\documentclass[reqno, 12pt]{amsart}
\usepackage{a4wide}
\usepackage{amscd}
\usepackage{enumerate}
\input xy
\xyoption{all}

\newtheorem*{thm A}{Theorem~A}
\newtheorem*{thm B}{Theorem~B}
\newtheorem*{thm C}{Theorem~C}
\newtheorem*{thm D}{Theorem~D}
\newtheorem{thm}{Theorem}[section]

\newtheorem{lm}[thm]{Lemma}
\newtheorem{re}[thm]{Remark}

\newtheorem*{MT}{Main Theorem}
\newtheorem*{MT1}{Theorem 1}
\newtheorem*{MT2}{Theorem 2}

\def \NBo{SU_{2,m-1}/S(U_2{\cdot}U_{m-1})}
\def \NBt{SU_{2,m}/S(U_{2}{\cdot}U_{m})}

\numberwithin{equation}{section}

\begin{document}

\title[Complex quadric with commuting Ricci tensor]{Real hypersurfaces in the complex quadric with commuting Ricci tensor}
\author{\textsc{Young Jin Suh and Doo Hyun Hwang}}

\address{Kyungpook National University \\ College of Natural Sciences\\  Department of Mathematics \\
Daegu 702-701\\ Republic of Korea}
\email{yjsuh@knu.ac.kr}
\email{engus0322@knu.ac.kr}
\date{}

\begin{abstract}
We introduce the notion of commuting Ricci tensor for real hypersurfaces in the complex quadric $Q^m = SO_{m+2}/SO_mSO_2$ . It is shown that the commuting Ricci tensor gives that the unit normal vector field $N$ becomes
$\frak A$-principal or $\frak A$-isotropic. Then according to each case, we give a complete classification of real hypersurfaces in $Q^m = SO_{m+2}/SO_mSO_2$ with commuting Ricci tensor.
\end{abstract}

\maketitle
\thispagestyle{empty}

\footnote[0]{
This work was supported by grant Proj. No. NRF-2015-R1A2A1A-01002459 from National Research Foundation of Korea.}
\par
\vskip 6pt
\section{Introduction}\label{section 1}

\par
\vskip 6pt
In a class of Hermitian symmetric spaces of rank 2, usually we can give examples of
 Riemannian symmetric spaces $SU_{m+2}/S(U_2U_m)$ and $SU_{2,m}/S(U_2U_m)$, which are said to be complex two-plane Grassmannians and complex hyperbolic two-plane Grassmannians respectively (see \cite{S1}, \cite{S2},  \cite{SHA}, \cite{SH} and \cite{S3}). These are viewed as Hermitian symmetric spaces and quaternionic K\"{a}hler symmetric spaces equipped with the K\"{a}hler structure $J$ and  the quaternionic K\"{a}hler structure ${\mathfrak J}$ and they have rank $2$.
\par
\vskip 6pt

 Among the other different type of Hermitian symmetric space with rank $2$ in the class of compact type, we can give the example of complex quadric $Q^m = SO_{m+2}/SO_mSO_2$, which is a complex hypersurface in complex projective space ${\Bbb C}P^{m+1}$ (see  Suh \cite{SRP}, \cite{S5}, and Smyth \cite{SM}). The complex quadric can also be regarded as a kind of real Grassmann manifolds of compact type with rank 2 (see Kobayashi and Nomizu \cite{KN}).  Accordingly, the complex quadric admits two important geometric structures, a complex conjugation structure $A$ and a K\"ahler structure $J$, which anti-commute with each other, that is, $AJ=-JA$. Then for $m{\ge}2$ the triple $(Q^m,J,g)$ is a Hermitian symmetric space of compact type with rank $2$ and its maximal sectional curvature is equal to $4$ (see Klein \cite{K} and Reckziegel \cite{R}).
 \par
 \vskip 6pt
In the complex projective space ${\mathbb C}P^{m+1}$ and the quaternionic projective space ${\mathbb Q}P^{m+1}$ some classifications related to commuting Ricci tensor or commuting structure Jacobi operator were investigated by Kimura \cite{KM1}, \cite{KM2}, P\'erez \cite{P} and P\'erez and Suh \cite{PS1}, \cite{PS2} respectively.  Under the invariance of the shape operator along some distributions a new classification in the complex $2$-plane Grassmannian $G_2({\mathbb C}^{m+2})= SU_{m+2}/S(U_mU_2)$ was investigated.  By using this classification P\'erez and Suh \cite{PS3} proved a non-existence property for Hopf hypersurfaces in $G_2({\mathbb C}^{m+2})$ with parallel and commuting Ricci tensor.  Recently, Hwang, Lee and Woo \cite{HLW} considered the notion of semi-parallel symmetric operators and obtained a complete classificationn for Hopf hypersurfaces in $G_2({\mathbb C}^{m+2})$. Moreover, Suh \cite{S1} strengthened this result to hypersurfaces in $G_2({\mathbb C}^{m+2})$ with commuting Ricci tensor and gave a characterization of real hypersurfaces in $G_2({\mathbb C}^{m+2})= SU_{m+2}/S(U_mU_2)$ as follows:
 \par
 \vskip 6pt
\begin{thm A}\label{Theorem A}
Let $M$ be a Hopf real hypersurface in $G_2({\mathbb C}^{m+2})$ with commuting Ricci tensor, $m{\ge}3$. Then $M$ is locally congruent to a tube of radius $r$ over a totally geodesic $G_2({\mathbb C}^{m+1})$
in $G_2({\mathbb C}^{m+2})$.
\end{thm A}
\par
 \vskip 6pt
 Moreover, Suh \cite{S3} studied another classification for Hopf hypersurfaces in complex hyperbolic two-plane Grassmannians $SU_{2,m}/S(U_2U_m)$ with commuting Ricci tensor as follows:
\par
 \vskip 6pt
\begin{thm B}\label{Theorem B}
Let $M$ be a Hopf hypersurface in $\NBt$ with commuting Ricci tensor, $m{\ge}3$. Then $M$ is locally congruent to an open part of a
tube around some totally geodesic $\NBo$ in $\NBt$ or a horosphere
whose center at infinity with $JX{\in}{\frak J}X$ is singular.
\end{thm B}
\par
\vskip 6pt
    It is known that the Reeb flow on a real hypersurface in $G_2({\mathbb C}^{m+2})$ is isometric if and only if $M$ is an open part of a tube around a totally geodesic $ G_2({\mathbb C}^{m+1})\subset G_2({\mathbb C}^{m+2})$. Correponding to this result, in \cite{SH} we asserted that the Reeb flow on a real hypersurface in $SU_{2,m}/S(U_2U_m)$ is isometric if and only if $M$ is an open part of a tube around a totally geodesic
    $SU_{2,m-1}/S(U_2U_{m-1})\subset SU_{2,m}/S(U_2U_m)$. Here, the Reeb flow on real hypersurfaces in $SU_{m+2}/S(U_mU_2)$ or $SU_{2,m}/S(U_2U_m)$ is said to be {\it isometric} if the shape operator commutes with the structure tensor. In papers \cite{SRP} and \cite{S5} due to Suh, we have introduced this problem for real hypersurfaces in the complex quadric $Q^m = SO_{m+2}/SO_mSO_2$ and obtained the following result:
 \par
 \vskip 6pt
 \begin{thm C}\label{Theorem C}
 Let $M$ be a real hypersurface of the complex quadric $Q^m$, $m\geq 3$. The Reeb flow on $M$ is isometric if and only if $m$ is even, say $m = 2k$, and $M$ is an open part of a tube around a totally geodesic ${\mathbb C}P^k \subset Q^{2k}$.
\end{thm C}
\par
\vskip 6pt
Now at each point $z \in M$ let us consider a maximal ${\mathfrak A}$-invariant subspace ${\mathcal Q}_z$ of $T_zM$, $z{\in}M$, defined by
$${\mathcal Q}_z = \{X \in T_zM \mid AX \in T_zM\ {\rm for\ all}\ A \in {\mathfrak A}_z\}$$
of $T_zM$, $z{\in}M$. Thus for a case where the unit normal vector field $N$ is $\frak A$-isotropic it can be easily checked that the
orthogonal complement ${\mathcal Q}_z^{\bot}={\mathcal C}_z{\ominus}{\mathcal Q}_z$, $z{\in}M$, of the distribution $\mathcal Q$ in the complex subbundle ${\mathcal C}$, becomes
${\mathcal Q}_z^{\bot}=\text{Span}\{A{\xi},AN\}$. Here it can be easily checked that the vector fields $A{\xi}$ and $AN$ belong to the tangent space $T_zM$, $z{\in}M$ if the unit normal vector field $N$ becomes
$\frak A$-isotropic. Thus for a case where the unit normal vector field $N$ is $\frak A$-isotropic it can be easily checked that the
orthogonal complement ${\mathcal Q}_z^{\bot}={\mathcal C}_z{\ominus}{\mathcal Q}_z$, $z{\in}M$, of the distribution $\mathcal Q$ in the complex subbundle ${\mathcal C}$, becomes
${\mathcal Q}_z^{\bot}=\text{Span}\{A{\xi},AN\}$. Moreover, the vector fields $A{\xi}$ and $AN$ belong to the tangent space $T_zM$, $z{\in}M$ if the unit normal vector field $N$ becomes
$\frak A$-isotropic. Then motivated by the above result, in \cite{S5} we gave another theorem for real hypersurfaces in the complex quadric $Q^m$ with parallel Ricci tensor and $\frak A$-isotropic unit normal.
\par
\vskip 6pt

Apart from the complex structure $J$
there is another distinguished geometric structure on ${Q}^m$, namely a parallel rank two vector bundle ${\mathfrak A}$ which contains an $S^1$-bundle of real structures, that is, complex conjugations $A$ on the tangent spaces of $Q^m$. This geometric structure determines a maximal ${\mathfrak A}$-invariant subbundle ${\mathcal Q}$ of the tangent bundle $TM$ of a real hypersurface $M$ in ${Q}^m$.
\par
\vskip 6pt
When we consider a hypersurface $M$ in the complex quadric ${Q}^m$, under the assumption of some geometric properties the unit normal vector field $N$ of $M$ in $Q^m$ can be divided into two classes if either $N$ is $\mathfrak A$-isotropic or $\frak A$-principal (see \cite{SRP}, and \cite{S5}). In the first case where $N$ is $\mathfrak A$-isotropic, it was known that $M$ is locally congruent to a tube over a totally geodesic ${\mathbb C}P^k$ in $Q^{2k}$.
In the second case, when the unit normal $N$ is $\mathfrak A$-principal, we proved that a contact hypersurface $M$ in $Q^m$ is locally congruent to a tube over a totally geodesic and totally real submanifold $S^m$ in $Q^m$ (see \cite{S5}).
\par
\vskip 6pt
In the study of complex two-plane Grassmannian $G_2({\Bbb C}^{m+2})$ or complex hyperbolic two-plane Grassmannian $\NBt$ we studied hypersurfaces with parallel Ricci tensor and gave non-existence properties respectively
(see \cite{S2} and \cite{SW}). In \cite{S5} we also considered the notion of parallel Ricci tensor ${\nabla}\text{Ric}=0$ for hypersurfaces $M$ in $Q^m$. As a generalization of such facts, we consider the notion of harmonic curvature, that is, $({\nabla}_X\text{Ric})Y=({\nabla}_Y\text{Ric})X$ for any tangent vector fields $X$ and $Y$ on $M$ in $Q^m$ and proved the following (see \cite{S6})

\par
\vskip 6pt
\begin{thm D}\label{Theorem D}
Let $M$ be a Hopf real hypersurface in the complex quadric $Q^m$, $m \ge 4$, with harmonic curvature and $\frak A$-isotropic unit normal $N$. If the shape operator commutes with the structure tensor on the distribution ${\mathcal Q}^{\bot}$, then $M$ is locally congruent to an open part of a tube around $k$-dimensional complex projective space ${\mathbb C}P^k$ in $Q^m$, $m=2k$, or $M$ has at most $6$ distinct constant principal curvatures given by
$${\alpha}, \quad {\gamma}=0({\alpha}), \quad {\lambda}_1, \quad {\mu}_1, \quad {\lambda}_2\,\quad \text{and}\quad \mu_2 $$
with corresponding principal curvature spaces
$$T_{\alpha}=[{\xi}], \quad  T_{\gamma}=[A{\xi}, AN], \quad {\phi}(T_{\lambda_1})=T_{\mu_1}, \quad {\phi}T_{\lambda_2}= T_{\mu_2}.$$
$$\text{dim}\ T_{\lambda_1} + \text{dim}T_{\lambda_2}=m-2, \quad \text{dim}\ T_{\mu_1} + \text{dim}T_{\mu_2}=m-2.$$
Here four roots $\lambda_i$ and $\mu_i$, $i=1,2$ satisfy the equation
$$2x^2-2{\beta}x+2+{\alpha}{\beta}=0,$$
where the function $\beta$ denotes ${\beta}=\frac{{\alpha}^2+2{\pm}{\sqrt{({\alpha}^2+2)^2+4{\alpha}h}}}{\alpha}$. In particular, ${\alpha}=\sqrt{\frac{2m-1}{2}}$, $\gamma(={\alpha})=\sqrt{\frac{2m-1}{2}}$, ${\lambda}=0$, ${\mu}=-\frac{2{\sqrt 2}}{\sqrt{2m-1}}$, with multiplicities $1$, $2$, $m-2$ and $m-2$ respectively.
\end{thm D}
\par
\vskip 6pt

But from the assumption of harmonic curvature, it was impossible to derive the fact that either the unit normal $N$ is $\frak A$-isotropic or $\frak A$-principal. So in \cite{S6} we gave a complete classification with the further assumption of $\frak A$-isotropic as in Theorem D. For the case where the unit normal vector field $N$ is $\frak A$-principal we have proved that real hypersurfaces in $Q^m$ with harmonic curvature can not be existed.
\par
\vskip 6pt

But fortunately when we consider Ricci commuting, that is,
$\text{Ric}{\cdot}{\phi}={\phi}{\cdot}\text{Ric}$ for hypersurfaces $M$ in $Q^m$, we can assert that the unit normal vector field $N$ becomes either $\frak A$-isotropic or $\frak A$-principal.
Then motivated by such a result and using Theorem C, in this paper we give a complete classification for real hypersurfaces in the complex quadric $Q^m$ with commuting Ricci tensor, that is, $\text{Ric}{\cdot}{\phi}={\phi}{\cdot}\text{Ric}$ as follows:
\par
\vskip 6pt
\begin{MT}\label{Main Theorem}\quad Let $M$ be a Hopf real hypersurface in the complex quadric $Q^m$, $m{\ge}4$, with commuting Ricci tensor. If the shape operator commutes with the structure tensor on the distribution ${\mathcal Q}^{\bot}$, then $M$ is locally congruent to an open part of a tube around totally geodesic ${\Bbb C}P^k$ in $Q^{2k}$, $m=2k$ or $M$ has $3$ distinct constant principal curvatures given by
$${\alpha}={\sqrt {2(m-3)}}, {\gamma}=0, {\lambda}=0,\  \text{and}\ \mu=-{\frac{2}{\sqrt {2(m-3)}}}\ \text{or}$$
$${\alpha}={\sqrt {\frac{2}{3}(m-3)}}, {\gamma}=0, {\lambda}=0,\  \text{and}\ \mu=- \frac{\sqrt 6}{\sqrt {m-3}}$$
with corresponding principal curvature spaces respectively
$$T_{\alpha}=[{\xi}], T_{\gamma}=[A{\xi}, AN], {\phi}(T_{\lambda})=T_{\mu}, \  \text{and}\  \text{dim}\ T_{\lambda}=\text{dim}\ T_{\mu}=m-2.$$.
\end{MT}
\par
\vskip 6pt

\begin{re}
In the main theorem the second and the third ones can be explained geometrically as follows: the real hypersurface $M$ is locally congruent to $M_1{\times}{\Bbb C}$, where $M_1$ is a tube of radius $r=\frac{1}{\sqrt 2}\tan^{-1}{\sqrt{m-3}}$ or respectively, of radius $r=\frac{1}{\sqrt 2}\tan^{-1}{\sqrt{\frac{m-3}{3}}}$, around $(m-1)$-dimensional sphere $S^{m-1}$ in $Q^{m-1}$. That is , $M_1$ is a contact hypersurface defined by
$S{\phi}+{\phi}S=k{\phi}$, $k=-\frac{2}{\sqrt{2(m-3)}}$, and $k=-\frac{\sqrt 6}{\sqrt{m-3}}$ respectively (see Suh \cite{S5}). By the Segre embedding, the embedding  $M_1{\times}{\Bbb C}{\subset}
Q^{m-1}{\times}{\Bbb C}{\subset}Q^m$ is defined by $(z_0, z_1, {\cdots}, z_m, w){\to}(z_0w, z_1w,{\cdots}, z_mw, 0)$. Here $(z_0w)^2+(z_1w)^2+{\cdots}+(z_mw)^2=(z_0^2+{\cdots}+z_m^2)w^2=0$, where $\{z_0,{\cdots},z_{m}\}$
denotes a coordinate system in $Q^{m-1}$ satisfying $z_0^2+{\cdots}+z_m^2=0$.
\end{re}
\par
\vskip 6pt

Our paper is organized as follows. In section 2 we present basic material about the complex quadric ${Q}^m$, including its Riemannian curvature tensor and a description of its singular vectors of $Q^m$ like $\mathfrak A$-principal or $\mathfrak A$-isotropic unit normal vector field.  In section 3, we investigate the geometry of this subbundle ${\mathcal Q}$ for
hypersurfaces in $Q^m$ and some equations including Codazzi and fundamental formulas related to the vector fields $\xi$, $N$, $A{\xi}$, and $AN$ for the complex conjugation $A$ of $M$ in $Q^m$.
\par
\vskip 6pt
 In section 4, the first step is to derive the formula of Ricci commuting from the equation of Gauss for real hypersurfaces $M$ in $Q^m$ and to get a key lemma that the unit normal vector field $N$ can be divided into two clssses of normal vector suct that $N$ is either $\frak A$-isotropic or $\frak A$-principal, and show that a real hypersurface in $Q^m$, $m=2k$, which is a tube over a totally geodesic ${\Bbb C}P^k$ in $Q^{2k}$ naturally admits a commuting Ricci tensor. In sections 5, by the expressions of the shape operator $S$ for real hypersurfaces $M$ in $Q^m$, we present the proof of Main Theorem with $\frak A$-isotropic unit normal.
\par
\vskip 6pt
In section 6, we give a complete proof of Main Theorem with $\frak A$-principal unit normal. The first part of this proof is devoted to give some fundamental formulas from Ricci commuting and $\mathfrak A$-principal unit normal vector field. Then in the latter part of
the proof we will use the decomposition of two eigenspaces of the complex conjugation $A$ in $Q^m$ such that $T_zM=V(A){\oplus}JV(A)$, where two eigenspaces are defined by $V(A)=\{X{\in}T_zQ^m{\vert}AX=X\}$ and $JV(A)=\{X{\in}T_zQ^m{\vert}AX=-X\}$ respectively.
\par
\vskip 6pt

\section{The complex quadric}\label{section 2}
\par
\vskip 6pt
For more background to this section we refer to \cite{K}, \cite{KN}, \cite{R}, \cite{SRP}, and \cite{S5}. The complex quadric $Q^m$ is the complex hypersurface in ${\mathbb C}P^{m+1}$ which is defined by the equation $z_0^2 + \cdots + z_{m+1}^2 = 0$, where $z_0,\ldots,z_{m+1}$ are homogeneous coordinates on ${\mathbb C}P^{m+1}$. We equip $Q^m$ with the Riemannian metric $g$ which is induced from the Fubini-Study metric $\bar g$ on ${\mathbb C}P^{m+1}$ with constant holomorphic sectional curvature $4$. The Fubini-Study metric $\bar g$ is defined by ${\bar g}(X,Y)={\Phi}(JX,Y)$ for any vector fields $X$ and $Y$ on ${\mathbb C}P^{m+1}$ and a globally closed $(1,1)$-form $\Phi$ given by ${\Phi}=-4i{\partial}{\bar{\partial}}{\text log}f_j$ on an open set $U_j=\{[z^0, z^1,{\cdots},z^{m+1}]{\in}{\mathbb C}P^{m+1}{\vert}z^j{\not =}0\}$, where the function $f_j$ denotes $f_j={\sum}_{k=0}^{m+1}t_j^k{\bar t}_j^k$, and  $t_j^k=\frac{z^k}{z^j}$ for $j,k=0,{\cdots},m+1$. Then naturally the K\"{a}hler structure on ${\mathbb C}P^{m+1}$ induces canonically a K\"{a}hler structure $(J,g)$ on the complex quadric $Q^m$.
\par
\vskip 6pt
The complex projective space ${\mathbb C}P^{m+1}$ is a Hermitian symmetric space of the special unitary group $SU_{m+2}$, namely ${\mathbb C}P^{m+1} = SU_{m+2}/S(U_{m+1}U_1)$. We denote by $o = [0,\ldots,0,1] \in {\mathbb C}P^{m+1}$ the fixed point of the action of the stabilizer $S(U_{m+1}U_1)$. The special orthogonal group $SO_{m+2} \subset SU_{m+2}$ acts on ${\mathbb C}P^{m+1}$ with cohomogeneity one. The orbit containing $o$ is a totally geodesic real projective space ${\mathbb R}P^{m+1} \subset {\mathbb C}P^{m+1}$. The second singular orbit of this action is the complex quadric $Q^m = SO_{m+2}/SO_mSO_2$. This homogeneous space model leads to the geometric interpretation of the complex quadric $Q^m$ as the Grassmann manifold $G_2^+({\mathbb R}^{m+2})$ of oriented $2$-planes in ${\mathbb R}^{m+2}$. It also gives a model of $Q^m$ as a Hermitian symmetric space of rank $2$. The complex quadric $Q^1$ is isometric to a sphere $S^2$ with constant curvature, and $Q^2$ is isometric to the Riemannian product of two $2$-spheres with constant curvature. For this reason we will assume $m \geq 3$ from now on.
\par
\vskip 6pt
For a nonzero vector $z \in {\mathbb C}^{m+2}$ we denote by $[z]$ the complex span of $z$, that is, $[z] = \{\lambda z \mid \lambda \in {\mathbb C}\} $.
Note that by definition $[z]$ is a point in ${\mathbb C}P^{m+1}$.
As usual, for each $[z] \in {\mathbb C}P^{m+1}$ we identify $T_{[z]}{\mathbb C}P^{m+1}$ with the orthogonal complement ${\mathbb C}^{m+2} \ominus [z]$ of $[z]$ in ${\mathbb C}^{m+2}$. For $[z] \in {Q}^m$ the tangent space $T_{[z]}{Q}^m$ can then be identified canonically with the orthogonal complement ${\mathbb C}^{m+2} \ominus ([z] \oplus [\bar{z}])$ of $[z] \oplus [\bar{z}]$ in ${\mathbb C}^{m+2}$ (see Kobayashi and Nomizu \cite{KN}). Note that $\bar{z} \in \nu_{[z]}{Q}^m$ is a unit normal vector of ${Q}^m$ in ${\mathbb C}P^{m+1}$ at the point $[z]$.
\par
\vskip 6pt
We denote by $A_{\bar{z}}$ the shape operator of ${Q}^m$ in ${\mathbb C}P^{m+1}$ with respect to the unit normal $\bar{z}$. It is defined by $A_{\bar{z}}w = {\bar{\nabla}}_w{\bar z}= {\bar w}$  for a complex Euclidean connection $\bar{\nabla}$ induced from ${\mathbb C}^{m+2}$ and all $w \in T_{[z]}Q^m$. That is, the shape operator $A_{\bar{z}}$ is just a complex conjugation restricted to $T_{[z]}{Q}^m$. Moreover, it satisfies the following for any $w \in T_{[z]}Q^m$ and any ${\lambda}{\in}S^1{\subset}{\mathbb C}$
\begin{equation*}
\begin{split}
A_{{\lambda}{\bar z}}^2w=&A_{{\lambda}{\bar z}}A_{{\lambda}{\bar z}}w=A_{{\lambda}{\bar z}}{\lambda}{\bar w}\\
=&{\lambda}A_{\bar z}{\lambda}{\bar w}={\lambda}{\bar{\nabla}}_{{\lambda}{\bar w}}{\bar z}={\lambda}{\bar{\lambda}}{\bar{\bar w}}\\
=&{\vert}{\lambda}{\vert}^2w=w.
\end{split}
\end{equation*}
Accordingly, $A_{{\lambda}{\bar z}}^2=I$ for any ${\lambda}{\in}S^1$. So the shape operator $A_{\bar{z}}$ becomes an anti-commuting involution such that $A_{\bar{z}}^2=I$ and $AJ=-JA$ on the complex vector space $T_{[z]}{Q}^m$ and
$$T_{[z]}{Q}^m = V(A_{\bar{z}}) \oplus JV(A_{\bar{z}}),$$
where $V(A_{\bar{z}}) =  {\mathbb R}^{m+2} \cap T_{[z]}{Q}^m$ is the $(+1)$-eigenspace and $JV(A_{\bar{z}}) = i{\mathbb R}^{m+2} \cap T_{[z]}Q^m$ is the $(-1)$-eigenspace of $A_{\bar{z}}$. That is, $A_{\bar{z}}X=X$ and $A_{\bar{z}}JX=-JX$, respectively, for any $X{\in}V(A_{\bar{z}})$.
\par
\vskip 6pt

Geometrically this means that the shape operator $A_{\bar{z}}$ defines a real structure on the complex vector space $T_{[z]}Q^m$, or equivalently, is a complex conjugation on $T_{[z]}Q^m$. Since the real codimension of $Q^m$ in ${\mathbb C}P^{m+1}$ is $2$, this induces an $S^1$-subbundle ${\mathfrak A}$ of the endomorphism bundle ${\rm End}(TQ^m)$ consisting of complex conjugations.
\par
\vskip 6pt
There is a geometric interpretation of these conjugations. The complex quadric $Q^m$ can be viewed as the complexification of the $m$-dimensional sphere $S^m$. Through each point $[z] \in Q^m$ there exists a one-parameter family of real forms of $Q^m$ which are isometric to the sphere $S^m$. These real forms are congruent to each other under action of the center $SO_2$ of the isotropy subgroup of $SO_{m+2}$ at $[z]$. The isometric reflection of $Q^m$ in such a real form $S^m$ is an isometry, and the differential at $[z]$ of such a reflection is a conjugation on $T_{[z]}Q^m$. In this way the family ${\mathfrak A}$ of conjugations on $T_{[z]}Q^m$ corresponds to the family of real forms $S^m$ of $Q^m$ containing $[z]$, and the subspaces $V(A) \subset T_{[z]}Q^m$ correspond to the tangent spaces $T_{[z]}S^m$ of the real forms $S^m$ of $Q^m$.
\par
\vskip 6pt
The Gauss equation for $Q^m \subset {\mathbb C}P^{m+1}$ implies that the Riemannian curvature tensor $\bar R$ of $Q^m$ can be described in terms of the complex structure $J$ and the complex conjugations $A \in {\mathfrak A}$:
\begin{eqnarray*}
{\bar R}(X,Y)Z & = & g(Y,Z)X - g(X,Z)Y + g(JY,Z)JX - g(JX,Z)JY - 2g(JX,Y)JZ \\
 & & + g(AY,Z)AX - g(AX,Z)AY + g(JAY,Z)JAX - g(JAX,Z)JAY.
\end{eqnarray*}
Note that $J$ and each complex conjugation $A$ anti-commute, that is, $AJ = -JA$ for each $A \in {\mathfrak A}$.
\par
\vskip 6pt
Recall that a nonzero tangent vector $W \in T_{[z]}Q^m$ is called singular if it is tangent to more than one maximal flat in $Q^m$. There are two types of singular tangent vectors for the complex quadric $Q^m$:
\begin{itemize}
\item[1.] If there exists a conjugation $A \in {\mathfrak A}$ such that $W \in V(A)$, then $W$ is singular. Such a singular tangent vector is called ${\mathfrak A}$-principal.
\item[2.] If there exist a conjugation $A \in {\mathfrak A}$ and orthonormal vectors $X,Y \in V(A)$ such that $W/||W|| = (X+JY)/\sqrt{2}$, then $W$ is singular. Such a singular tangent vector is called ${\mathfrak A}$-isotropic.
\end{itemize}
\par
\vskip 6pt
For every unit tangent vector $W \in T_{[z]}Q^m$ there exist a conjugation $A \in {\mathfrak A}$ and orthonormal vectors $X,Y \in V(A)$ such that
\[
W = \cos(t)X + \sin(t)JY
\]
for some $t \in [0,\pi/4]$. The singular tangent vectors correspond to the values $t = 0$ and $t = \pi/4$. If $0 < t < \pi/4$ then the unique maximal flat containing $W$ is ${\mathbb R}X \oplus {\mathbb R}JY$.
Later we will need the eigenvalues and eigenspaces of the Jacobi operator $R_W = R(\cdot,W)W$ for a singular unit tangent vector $W$.
\par
\vskip 6pt
\begin{itemize}
\item[1.] If $W$ is an ${\mathfrak A}$-principal singular unit tangent vector with respect to $A \in {\mathfrak A}$, then the eigenvalues of $R_W$ are $0$ and $2$ and the corresponding eigenspaces are ${\mathbb R}W \oplus J(V(A) \ominus {\mathbb R}W)$ and $(V(A) \ominus {\mathbb R}W) \oplus {\mathbb R}JW$, respectively.
\item[2.] If $W$ is an ${\mathfrak A}$-isotropic singular unit tangent vector with respect to $A \in {\mathfrak A}$ and $X,Y \in V(A)$, then the eigenvalues of $R_W$ are $0$, $1$ and $4$ and the corresponding eigenspaces are ${\mathbb R}W \oplus {\mathbb C}(JX+Y)$, $T_{[z]}Q^m \ominus ({\mathbb C}X \oplus {\mathbb C}Y)$ and ${\mathbb R}JW$, respectively.
\end{itemize}

\par
\vskip 6pt

\section{Some general equations}\label{section 3}
\par
\vskip 6pt
Let $M$ be a  real hypersurface in $Q^m$ and denote by $(\phi,\xi,\eta,g)$ the induced almost contact metric structure. Note that $\xi = -JN$, where $N$ is a (local) unit normal vector field of $M$ and $\eta$ the corresponding $1$-form defined by ${\eta}(X)=g({\xi},X)$ for any tangent vector field $X$ on $M$. The tangent bundle $TM$ of $M$ splits orthogonally into  $TM = {\mathcal C} \oplus {\mathbb R}\xi$, where ${\mathcal C} = {\rm ker}(\eta)$ is the maximal complex subbundle of $TM$. The structure tensor field $\phi$ restricted to ${\mathcal C}$ coincides with the complex structure $J$ restricted to ${\mathcal C}$, and $\phi \xi = 0$.
\par
\vskip 6pt
At each point $z \in M$ we define a maximal ${\mathfrak A}$-invariant subspace of $T_zM$, $z{\in}M$ as follows:

\[
{\mathcal Q}_z = \{X \in T_zM \mid AX \in T_zM\ {\rm for\ all}\ A \in {\mathfrak A}_z\}.
\]

Then we want to introduce an important lemma which will be used in the proof of our main Theorem in the introduction.
\par
\vskip 6pt

\begin{lm}(\cite{SRP} and \cite{S5})\label{lemma 3.1}
For each $z \in M$ we have
\begin{itemize}
\item[(i)] If $N_z$ is ${\mathfrak A}$-principal, then ${\mathcal Q}_z = {\mathcal C}_z$.
\item[(ii)] IF $N_z$ is not ${\mathfrak A}$-principal, there exist a conjugation $A \in {\mathfrak A}$ and orthonormal vectors $X,Y \in V(A)$ such that $N_z = \cos(t)X + \sin(t)JY$ for some $t \in (0,\pi/4]$.
Then we have ${\mathcal Q}_z = {\mathcal C}_z \ominus {\mathbb C}(JX + Y)$.
\end{itemize}
\end{lm}
\par
\vskip 6pt
We now assume that $M$ is a Hopf hypersurface. Then we have
\[
S\xi = \alpha \xi ,
\]
where $S$ denotes the shape operator of the real hypersurfaces $M$ with the smooth function $\alpha = g(S\xi,\xi)$ on $M$.
When we consider the transform $JX$ by the K\"{a}hler structure $J$ on $Q^m$ for any vector field $X$ on $M$ in $Q^m$, we may put
$$JX={\phi}X+{\eta}(X)N$$
for a unit normal $N$ to $M$. Then we now consider the Codazzi equation
\begin{equation}\label{e31}
\begin{split}
g((\nabla_XS)Y - (\nabla_YS)X,Z) & =  \eta(X)g(\phi Y,Z) - \eta(Y) g(\phi X,Z) - 2\eta(Z) g(\phi X,Y) \\
& \quad \ \   + g(X,AN)g(AY,Z) - g(Y,AN)g(AX,Z)\\
& \quad \ \   + g(X,A\xi)g(J AY,Z) - g(Y,A\xi)g(JAX,Z).
\end{split}
\end{equation}
Putting $Z = \xi$ in (\ref{e31}) we get
\begin{equation*}
\begin{split}
g((\nabla_XS)Y - (\nabla_YS)X,\xi) & =   - 2 g(\phi X,Y) \\
& \quad \ \  + g(X,AN)g(Y,A\xi) - g(Y,AN)g(X,A\xi)\\
& \quad \ \  - g(X,A\xi)g(JY,A\xi) + g(Y,A\xi)g(JX,A\xi).
\end{split}
\end{equation*}
On the other hand, we have
\begin{eqnarray*}
 & & g((\nabla_XS)Y - (\nabla_YS)X,\xi) \\
& = & g((\nabla_XS)\xi,Y) - g((\nabla_YS)\xi,X) \\
& = & (X\alpha)\eta(Y) - (Y\alpha)\eta(X) + \alpha g((S\phi + \phi
S)X,Y) - 2g(S \phi SX,Y).
\end{eqnarray*}
Comparing the previous two equations and putting $X = \xi$ yields
$$
Y\alpha  =  (\xi \alpha)\eta(Y)  - 2g(\xi,AN)g(Y,A\xi) +
2g(Y,AN)g(\xi,A\xi).
$$
Reinserting this into the previous equation yields
\begin{eqnarray*}
 & & g((\nabla_XS)Y - (\nabla_YS)X,\xi) \\
& = &  - 2g(\xi,AN)g(X,A\xi)\eta(Y) + 2g(X,AN)g(\xi,A\xi)\eta(Y) \\
& &  + 2g(\xi,AN)g(Y,A\xi)\eta(X) - 2g(Y,AN)g(\xi,A\xi)\eta(X) \\
& & + \alpha g((\phi S + S\phi)X,Y) - 2g(S \phi SX,Y) .
\end{eqnarray*}
Altogether this implies
\begin{equation}\label{e32}
\begin{split}
0 = & 2g(S \phi SX,Y) - \alpha g((\phi S + S\phi)X,Y) - 2 g(\phi X,Y) \\
 & + g(X,AN)g(Y,A\xi) - g(Y,AN)g(X,A\xi)\\
 & - g(X,A\xi)g(JY,A\xi) + g(Y,A\xi)g(JX,A\xi)\\
 & + 2g(\xi,AN)g(X,A\xi)\eta(Y) - 2g(X,AN)g(\xi,A\xi)\eta(Y) \\
 &  - 2g(\xi,AN)g(Y,A\xi)\eta(X) + 2g(Y,AN)g(\xi,A\xi)\eta(X).
\end{split}
\end{equation}
At each point $z \in M$ we can choose $A \in {\mathfrak A}_z$ such that
\[ N = \cos(t)Z_1 + \sin(t)JZ_2 \]
for some orthonormal vectors $Z_1,Z_2 \in V(A)$ and $0 \leq t \leq \frac{\pi}{4}$ (see Proposition 3 in \cite{R}). Note that $t$ is a function on $M$.
First of all, since $\xi = -JN$, we have
\begin{equation}\label{e33}
\begin{split}
N  = & \cos(t)Z_1 + \sin(t)JZ_2, \\
AN  = & \cos(t)Z_1 - \sin(t)JZ_2, \\
\xi = & \sin(t)Z_2 - \cos(t)JZ_1, \\
A\xi = & \sin(t)Z_2 + \cos(t)JZ_1.
\end{split}
\end{equation}
This implies $g(\xi,AN) = 0$ and hence
\begin{equation}\label{e34}
\begin{split}
0  = & 2g(S \phi SX,Y) - \alpha g((\phi S + S\phi)X,Y) - 2 g(\phi X,Y) \\
 & + g(X,AN)g(Y,A\xi) - g(Y,AN)g(X,A\xi)\\
 & - g(X,A\xi)g(JY,A\xi) + g(Y,A\xi)g(JX,A\xi)\\
 &  - 2g(X,AN)g(\xi,A\xi)\eta(Y) + 2g(Y,AN)g(\xi,A\xi)\eta(X).
\end{split}
\end{equation}

\par
\vskip 6pt

\section{Ricci commuting and a Key Lemma}\label{section 4}
\par
\vskip 6pt
By the equation of Gauss, the curvature tensor $R(X,Y)Z$ for a real hypersurface $M$ in $Q^m$ induced from the curvature tensor
$\bar R$ of $Q^m$ can be described in terms of the complex structure $J$ and the complex conjugation $A \in {\mathfrak A}$ as follows:
\begin{eqnarray*}
 R(X,Y)Z & = & g(Y,Z)X - g(X,Z)Y + g({\phi}Y,Z){\phi}X - g({\phi}X,Z){\phi}Y - 2g({\phi}X,Y){\phi}Z \\
 & & + g(AY,Z)AX - g(AX,Z)AY + g(JAY,Z)JAX - g(JAX,Z)JAY\\
 & & + g(SY,Z)SX-g(SX,Z)SY
\end{eqnarray*}
for any $X,Y,Z{\in}T_zM$, $z{\in}M$.

Now let us put
$$AX=BX+{\rho}(X)N,$$
for any vector field $X{\in}T_zQ^m$, $z{\in}M$, ${\rho}(X)=g(AX,N)$, where $BX$ and ${\rho}(X)N$ respectively denote the tangential and normal component of the vector field $AX$. Then $A{\xi}=B{\xi}+{\rho}({\xi})N$ and
${\rho}({\xi})=g(A{\xi},N)=0$. Then it follows that
\begin{equation*}
\begin{split}
AN=&AJ{\xi}=JA{\xi}=-J(B{\xi}+{\rho}({\xi})N)\\
=&-({\phi}B{\xi}+{\eta}(B{\xi})N).
\end{split}
\end{equation*}
The equation gives $g(AN,N)=-{\eta}(B{\xi})$ and $g(AN,{\xi})=0$. From this, together with the definition of the Ricci tensor, we have
\begin{equation*}
\begin{split}
Ric(X)=&(2m-1)X-3{\eta}(X){\xi}-g(AN,N)AX+g(AX,N)AN\\
&+{\eta}(AX)A{\xi}+(Tr S)SX-S^2X.
\end{split}
\end{equation*}
Then, summing up with the above formulas, it can be rearranged as follows:
\begin{equation*}
\begin{split}
Ric(X)=&(2m-1)X-3{\eta}(X){\xi}+{\eta}(B{\xi})\{BX+{\rho}(X)N\}\\
&+{\rho}(X)\{-{\phi}B{\xi}-{\eta}(B{\xi})N\}+{\eta}(BX)B{\xi}+(Tr S)SX-S^2X
\end{split}
\end{equation*}
From this, together with the assumption of Ricci commuting, that is, ${\phi}{\cdot}Ric(X)=Ric{\cdot}{\phi}X$, it follows that
\begin{equation}\label{e41}
\begin{split}
(2m-1){\phi}X+{\eta}(B{\xi}){\phi}BX-{\rho}(X){\phi}^2B{\xi}\\
+{\eta}(BX){\phi}B{\xi}+(Tr S){\phi}SX-{\phi}S^2X\\
=(2m-1){\phi}X+{\eta}(B{\xi})B{\phi}X-{\rho}({\phi}X){\phi}B{\xi}\\
+{\eta}(B{\phi}X)B{\xi}+(Tr S)S{\phi}X-S^2{\phi}X.
\end{split}
\end{equation}
Here  we want to use the following formulas
$$
{\eta}(BX)=g(A{\xi},X),$$
\begin{equation*}
\begin{split}
{\eta}(B{\phi}X)=&g(A{\xi},{\phi}X)=g(A{\xi},JX-{\eta}(X)N)=g(AJ{\xi},X)\\
&=g(AN,X)={\rho}(X),
\end{split}
\end{equation*}
\begin{equation*}
\begin{split}
{\rho}({\phi}X)=&g(A{\phi}X,N)=g(AJ{\phi}X,{\xi})=g(J{\phi}X,A{\xi})\\
=&g({\phi}^2X+{\eta}({\phi}X)N, A{\xi})=-g(X,A{\xi})+{\eta}(X)g({\xi},A{\xi}),
\end{split}
\end{equation*}
$$
{\rho}(X)={\eta}(B{\phi}X).$$

Summing up these formulas into (4.1), we have
\begin{equation}\label{e42}
\begin{split}
{\eta}&(B{\xi}){\phi}BX-{\eta}(B{\phi}X){\eta}(B{\xi}){\xi}+(Tr S){\phi}SX-{\phi}S^2X\\
=&{\eta}(B{\xi})B{\phi}X-{\eta}(X){\eta}(B{\xi}){\phi}B{\xi}+(Tr S)S{\phi}X-S^2{\phi}X.
\end{split}
\end{equation}
Then, by taking inner product of (4.2) with $\xi$ and using that $M$ is Hopf, it follows that
\begin{equation}\label{e43}
\begin{split}
{\eta}(B{\xi}){\phi}B{\xi}=0.
\end{split}
\end{equation}
Then the formula (\ref{e42}) becomes
\begin{equation}\label{e44}
{\eta}(B{\xi})({\phi}B-B{\phi})X+(Tr S)({\phi}S-S{\phi})X-({\phi}S^2-S^2{\phi})X=0.
\end{equation}
Here we want to give a remark as follows:
\par
\vskip 6pt
\begin{re}
Let $M$ be a real hypersurface over a totally geodesic ${\mathbb C}P^k{\subset}Q^{2k}$, $m=2k$. Then in papers due to \cite{SRP} and \cite{S5}
the structure tensor commutes with the shape operator, that is, $S{\phi}={\phi}S$. Moreover, the unit normal vector field $N$ becomes $\mathfrak A$-isotropic.
This gives ${\eta}(B{\xi})=g(A{\xi},{\xi})=0$. So it naturally satisfies the formula (\ref{e42}), that is, Ricci commuting.
\end{re}
\par
\vskip 6pt
On the other hand, from (\ref{e43}) we assert an important lemma as follows:
\par
\vskip 6pt
\begin{lm}\label{Lemma 4.2}\quad Let $M$ be a real hypersurface in $Q^m$, $m{\ge}3$, with commuting Ricci tensor.  Then the unit normal vector field $N$ becomes singular,
that is, $N$ is $\mathfrak A$-isotropic or $\mathfrak A$-principal.
\end{lm}
\par
\vskip 6pt
\begin{proof}\quad From (\ref{e43}) we get
$${\eta}(B{\xi})=0\quad \text{or}\quad {\phi}B{\xi}=0.$$
The first case gives that ${\eta}(B{\xi})=g(A{\xi},{\xi})=\cos 2t=0$, that is, $t=\frac{\pi}{4}$. This implies that the unit normal $N$ becomes
$N=\frac{X+JY}{\sqrt 2}$, which means that $N$ is $\mathfrak A$-isotropic.

The second case gives that
$${\rho}(X)=g(AX,N)={\eta}(B{\phi}X)=-g(X,{\phi}B{\xi})=0,$$
which means that $AX{\in}T_zM$ for any $A{\in}\mathfrak A$, $X{\in}T_zM$, $z{\in}M$.
This implies ${\mathcal Q}_z = {\mathcal C}_z$, $z{\in}M$, and $N$ is $\mathfrak A$-principal, that is, $AN=N$.
\end{proof}
\vskip 6pt

In order to prove our main theorem in the introduction, by virtue of Lemma 4.2, we can divide into two classes of hypersurfaces in $Q^m$ with the unit normal $N$ is $\mathfrak A$-principal or $\mathfrak A$-isotropic.
When $M$ is with $\mathfrak A$-isotropic, in section 5 we will give its proof in detail and in section 6 we will give the remainder proof for the case that $M$ has a $\mathfrak A$-principal normal vector field.

\par
\vskip 6pt

\section{Proof of Main Theorem with $\frak A$-isotropic}\label{section 5}
\par
\vskip 6pt
In this section we want to prove our Main Theorem for real hypersurfaces $M$ in $Q^m$ with commuting Ricci tensor when the unit normal vector field becomes $\frak A$-isotropic.
\par
\vskip 6pt
Since we assumed that the unit normal $N$ is $\frak A$-isotropic, by the definition in section 3 we know that $t=\frac{\pi}{4}$. Then by the expression
of the $\frak A$-isotropic unit normal vector field, (\ref{e33}) gives $N=\frac{1}{\sqrt 2}Z_1+\frac{1}{\sqrt 2}JZ_2$ . This implies that $g(A{\xi},{\xi})=0$. Since the unit normal $N$ is $\frak A$-isotropic, we know that $g({\xi},A{\xi})=0$. Moreover, by (\ref{e34}) and using an anti-commuting property $AJ=-JA$ between the complex conjugation $A$ and the K\"ahler structure $J$ , we proved the following (see also Lemma 4.2 in \cite{SRP})
\par
\vskip 6pt
\begin{lm}\label{lemma 5.1}
Let $M$ be a Hopf hypersurface in $Q^m$ with (local) unit normal vector field $N$. For each pont in $z \in M$ we choose $A \in {\mathfrak A}_z$ such that
$N_z = \cos(t)Z_1 + \sin(t)JZ_2$ holds
for some orthonormal vectors $Z_1,Z_2 \in V(A)$ and $0 \leq t \leq \frac{\pi}{4}$. Then
\begin{eqnarray*}
0 & = & 2g(S \phi SX,Y) - \alpha g((\phi S + S\phi)X,Y) -  2g(\phi X,Y) \\
& & + 2g(X,AN)g(Y,A\xi) - 2g(Y,AN)g(X,A\xi) \\
& &  + 2g(\xi,A\xi) \{g(Y,AN)\eta(X) - g(X,AN)\eta(Y) \}
\end{eqnarray*}
holds for all vector fields $X,Y$ on $M$.
\end{lm}
Then by virtue of $\frak A$-isotropic unit normal, Lemma 5.1 becomes the following
\begin{equation}\label{e51}
2S{\phi}SX={\alpha}(S{\phi}+{\phi}S)X+2{\phi}X-2g(X,AN)A{\xi}+2g(X,A{\xi})AN.
\end{equation}

Now let us consider the distribution ${\mathcal Q}^{\bot}$, which is an orthogonal complement of the maximal $\frak A$-invariant subspace $\mathcal Q$ in the complex subbundle $\mathcal C$ of $T_zM$, $z{\in}M$ in $Q^m$.
Then by Lemma \ref{lemma 3.1} in section 3, the orthogonal complement ${\mathcal Q}^{\bot}={\mathcal C}{\ominus}{\mathcal Q}$ becomes ${\mathcal C}{\ominus}{\mathcal Q}=\text{Span}\ [AN,A{\xi}]$.
From the assumption of $S{\phi}={\phi}S$ on the distribution ${\mathcal Q}^{\bot}$ it can be easily checked that the distribution ${\mathcal Q}^{\bot}$ is invariant by the shape operator $S$. Then $(\ref{e51})$ gives the following for $SAN={\lambda}AN$

\begin{eqnarray*}
(2{\lambda}-{\alpha})S{\phi}AN&=&({\alpha}{\lambda}+2){\phi}AN-2A{\xi}\\
&=&({\alpha}{\lambda}+2){\phi}AN-2{\phi}AN\\
&=&{\alpha}{\lambda}{\phi}AN.
\end{eqnarray*}
Then $A{\xi}={\phi}AN$ gives the following

\begin{equation}\label{e52}
SA{\xi}=\frac{{\alpha}{\lambda}}{2{\lambda}-{\alpha}}A{\xi}.
\end{equation}

Then from the assumption $S{\phi}={\phi}S$ on ${\mathcal Q}^{\bot}={\mathcal C}{\ominus}{\mathcal Q}$ it follows that ${\lambda}=\frac{{\alpha}{\lambda}}{2{\lambda}-{\alpha}}$ gives
\begin{equation}\label{e53}
{\lambda}=0\ \text{or}\  {\lambda}={\alpha}.
\end{equation}
\par
\vskip 6pt
On the other hand, on the distribution $\mathcal Q$ we know that $AX{\in}T_zM$, $z{\in}M$, because $AN{\in}Q$. So (\ref{e51}), together with the fact that $g(X,A{\xi})=0$ and $g(X,AN)=0$ for any $X{\in}{\mathcal Q}$, imply that
\begin{equation}\label{e54}
2S{\phi}SX = {\alpha}(S{\phi}+{\phi}S)X+2{\phi}X.
\end{equation}
Then we can take an orthonormal basis $X_1,{\cdots},X_{2(m-2)}{\in}{\mathcal Q}$ such that $AX_i={\lambda}_iX_i$ for $i=1,{\cdots},m-2$. Then by (\ref{e51}) we know that
$$S{\phi}X_i=\frac{{\alpha}{\lambda}_i+2}{2{\lambda_i}-{\alpha}}{\phi}X_i.$$
Accordingly, by (\ref{e53}) the shape operator $S$ can be expressed in such a way that

\begin{equation*}
S=\begin{bmatrix}
{\alpha} & 0 & 0 & 0 & {\cdots} & 0 & 0 & {\cdots} & 0\\
0 & 0({\alpha}) & 0 & 0 & {\cdots} & 0 & 0 & {\cdots} & 0 \\
0 & 0 & 0({\alpha}) & 0 & {\cdots} & 0 & 0 & {\cdots} & 0 \\
0 & 0 & 0 & {\lambda}_1 & {\cdots} & 0 & 0 & {\cdots} & 0\\
{\vdots} & {\vdots} & {\vdots} & {\vdots} & {\ddots} & {\vdots} & {\vdots} & {\cdots} & {\vdots}\\
0 & 0 & 0 & 0 & {\cdots} & {\lambda}_{m-2} & 0 & {\cdots} & 0\\
0 & 0 & 0 & 0 & {\cdots} & 0 & {\mu}_{1} & {\cdots} & 0 \\
{\vdots} & {\vdots} & {\vdots} & {\vdots} & {\vdots} & {\vdots} & {\vdots} & {\ddots} & {\vdots}\\
0 & 0 & 0 & 0 & {\cdots} & 0 & 0 & {\cdots} & {\mu}_{m-2}
\end{bmatrix}
\end{equation*}

From the equation (4.4), together with ${\eta}(B{\xi})=g(A{\xi},{\xi})=0$, we have that
\begin{equation}\label{e55}
h({\phi}S-S{\phi})X=({\phi}S^2-S^2{\phi})X,
\end{equation}
where $h= Tr S$ denotes the trace of the shape operator of $M$ in $Q^m$.

Now let us consider the Ricci commuting with $\mathfrak A$-isotropic normal for $SX={\lambda}X$, $X{\in}{\mathcal C}$. Then from (5.5) it follows that for $SX={\lambda}X$, $X{\in}{\mathcal C}$
\begin{equation}\label{e56}
({\lambda}-{\mu})\{h-({\lambda}+{\mu})\}{\phi}X=0,
\end{equation}
where we have used $S{\phi}X={\mu}{\phi}X$ for ${\mu}=\frac{{\alpha}{\lambda}+2}{2{\lambda}-{\alpha}}$, $X{\in}{\mathcal Q}$ and ${\mu}=\frac{{\alpha}{\lambda}}{2{\lambda}-{\alpha}}$, $X{\in}Q^{\bot}$ respectively.
Then (\ref{e56}) gives that
\begin{equation}\label{e57}
{\lambda}={\mu}\quad \text{or}\quad h={\lambda}+{\mu}.
\end{equation}

On the other hand, we consider the following for $SX={\lambda}X$, $X{\in}{\mathcal C}$. Then (\ref{e55}) gives
\begin{equation}\label{e58}
h{\lambda}{\phi}X-hS{\phi}X={\lambda}^2{\phi}X-S^2{\phi}X.
\end{equation}
Here we decompose $X{\in}{\mathcal C}={\mathcal Q}{\oplus}{\mathcal Q}^{\bot}$ in such a way that
$$X=Y+Z,$$
where $Y{\in}{\mathcal Q}$ and $Z{\in}{\mathcal Q}^{\bot}$. Then $SX={\lambda}X={\lambda}Y+{\lambda}Z$ gives the following
$$
SY={\lambda}Y\quad \text{and}\quad SZ={\lambda}Z,$$
because the distribution $\mathcal Q$ and ${\mathcal Q}^{\bot}$ are invariant by the shape operator. Then by using the matrix representation of the shape operator
the formula (\ref{e58}) gives the following decomposition
\begin{equation}\label{e59}
h{\lambda}{\phi}Y-h(\frac{{\alpha}{\lambda}+2}{2{\lambda}-{\alpha}}){\phi}Y={\lambda}^2{\phi}Y-(\frac{{\alpha}{\lambda}+2}{2{\lambda}-{\alpha}})^2{\phi}Y, \quad Y{\in}{\mathcal Q},
\end{equation}

\begin{equation}\label{e510}
h{\lambda}{\phi}Z-h(\frac{{\alpha}{\lambda}}{2{\lambda}-{\alpha}}){\phi}Z={\lambda}^2{\phi}Z-(\frac{{\alpha}{\lambda}}{2{\lambda}-{\alpha}})^2{\phi}Z, \quad Z{\in}{\mathcal Q}^{\bot}.
\end{equation}

By taking inner products of (\ref{e59}) and (\ref{e510}) with the vector fields ${\phi}Y$ and ${\phi}Z$ respectively, we have
\begin{equation}\label{e511}
{\lambda}^2-h{\lambda}+\frac{{\alpha}{\lambda}+2}{2{\lambda}-{\alpha}}\{h-\frac{{\alpha}{\lambda}+2}{2{\lambda}-{\alpha}}\}=0,
\end{equation}
\begin{equation}\label{e512}
{\lambda}^2-h{\lambda}+\frac{{\alpha}{\lambda}}{2{\lambda}-{\alpha}}\{h-\frac{{\alpha}{\lambda}}{2{\lambda}-{\alpha}}\}=0.
\end{equation}
Then substracting (\ref{e512}) from (\ref{e511}) gives
\begin{equation}\label{e513}
h=\frac{2{\alpha}{\lambda}+2}{2{\lambda}-{\alpha}}.
\end{equation}

Now from (\ref{e57}) let us consider the following two cases:
\par
\vskip 6pt
Case I. \quad ${\lambda}={\mu}$
\par
\vskip 6pt
From the matrix representation of the shape operator, ${\lambda}=\frac{{\alpha}{\lambda}+2}{2{\lambda}-{\alpha}}$ gives that
$${\lambda}^2-{\alpha}{\lambda}-1=0.$$
Since the discreminant $D={\alpha}^2+4>0$, we have two distinct solutions ${\lambda}=\cot r$ and ${\mu}=-\tan r$ with the mutiplicities $(m-2)$ and $(m-2)$ respectively.

That is, the shape operator $S$ can be expressed in such a way that

\begin{equation*}
S=\begin{bmatrix}
{\alpha} & 0 & 0 & 0 & {\cdots} & 0 & 0 & {\cdots} & 0\\
0 & 0({\alpha}) & 0 & 0 & {\cdots} & 0 & 0 & {\cdots} & 0 \\
0 & 0 & 0({\alpha}) & 0 & {\cdots} & 0 & 0 & {\cdots} & 0 \\
0 & 0 & 0 & {\cot r} & {\cdots} & 0 & 0 & {\cdots} & 0\\
{\vdots} & {\vdots} & {\vdots} & {\vdots} & {\ddots} & {\vdots} & {\vdots} & {\cdots} & {\vdots}\\
0 & 0 & 0 & 0 & {\cdots} & {\cot r} & 0 & {\cdots} & 0\\
0 & 0 & 0 & 0 & {\cdots} & 0 & -{\tan r} & {\cdots} & 0 \\
{\vdots} & {\vdots} & {\vdots} & {\vdots} & {\vdots} & {\vdots} & {\vdots} & {\ddots} & {\vdots}\\
0 & 0 & 0 & 0 & {\cdots} & 0 & 0 & {\cdots} & -{\tan r}
\end{bmatrix}
\end{equation*}

This means that the shape operator $S$ commutes with the structure tensor $\phi$, that is, $S{\cdot}{\phi}={\phi}{\cdot}S$. Then by Theorem C, $m=2k$, and $M$ is locally congruent to an open part of a tube around a totally geodesic
${\Bbb C}P^k$ in $Q^{2k}$.
\par
\vskip 6pt
Case II. \quad ${\lambda}{\not =}{\mu}$
\par
\vskip 6pt
Now we only consider ${\lambda}{\not =}{\mu}$ on the distribution $\mathcal Q$. Since on the distribution ${\mathcal Q}^{\bot}$ we have assumed that $S{\phi}={\phi}S$, so it follows that
${\lambda}=\frac{{\alpha}{\lambda}}{2{\lambda}-{\alpha}}$. This gives ${\lambda}=0$ or ${\lambda}={\alpha}$ on the distribution ${\mathcal Q}^{\bot}$. Moreover, by the Ricci commuting, we have
the following from (\ref{e53}), together with (\ref{e56}),(\ref{e57}) and (\ref{e513})

\begin{equation}
\begin{split}
h=&{\lambda}+{\mu}={\lambda}+\frac{{\alpha}{\lambda}+2}{2{\lambda}-{\alpha}}\\
=&\frac{2{\alpha}{\lambda}+2}{2{\lambda}-{\alpha}}.
\end{split}
\end{equation}

This gives ${\lambda}=0$ or ${\lambda}={\alpha}$. Then we can divide into two subcases as follows:
\par
\vskip 6pt
Subcase 2.1. ${\lambda}=0$.
\par
\vskip 6pt
Then we can arrange the matrix of the shape opeartor such that
\begin{equation*}
S=\begin{bmatrix}
{\alpha} & 0 & 0 & 0 & {\cdots} & 0 & 0 & {\cdots} & 0\\
0 & 0({\alpha}) & 0 & 0 & {\cdots} & 0 & 0 & {\cdots} & 0 \\
0 & 0 & 0({\alpha}) & 0 & {\cdots} & 0 & 0 & {\cdots} & 0 \\
0 & 0 & 0 & 0 & {\cdots} & 0 & 0 & {\cdots} & 0\\
{\vdots} & {\vdots} & {\vdots} & {\vdots} & {\ddots} & {\vdots} & {\vdots} & {\cdots} & {\vdots}\\
0 & 0 & 0 & 0 & {\cdots} & 0 & 0 & {\cdots} & 0\\
0 & 0 & 0 & 0 & {\cdots} & 0 & {-\frac{2}{\alpha}} & {\cdots} & 0 \\
{\vdots} & {\vdots} & {\vdots} & {\vdots} & {\vdots} & {\vdots} & {\vdots} & {\ddots} & {\vdots}\\
0 & 0 & 0 & 0 & {\cdots} & 0 & 0 & {\cdots} & {-\frac{2}{\alpha}}
\end{bmatrix}
\end{equation*}

In this case the formula $h={\lambda}+{\mu}$ and the notion of trace $h$ of the shape operator $S$ gives
\begin{equation*}
\begin{split}
h=&0-\frac{2}{\alpha}\\
=&{\alpha}+(m-2)(-\frac{2}{\alpha}).
\end{split}
\end{equation*}
Then it gives that ${\alpha}^2=2(m-3)$, that is, ${\alpha}={\sqrt {2(2m-3)}}$.

Now let us consider another case that $h=3{\alpha}+(m-2)(-\frac{2}{\alpha})$. Then, from this, together with $h={\lambda}+{\mu}=-\frac{2}{\alpha}$,
we know ${\alpha}^2=\frac{2}{3}(m-3)$. Then ${\alpha}=\sqrt{\frac{2}{3}(m-3)}$.
\par
\vskip 6pt

Subcase 2.2. ${\lambda}={\alpha}$.
\par
\vskip 6pt
In this subcase, the expression of the shape operator beomes

\begin{equation*}
S=\begin{bmatrix}
{\alpha} & 0 & 0 & 0 & {\cdots} & 0 & 0 & {\cdots} & 0\\
0 & {\alpha}(0) & 0 & 0 & {\cdots} & 0 & 0 & {\cdots} & 0 \\
0 & 0 & {\alpha}(0) & 0 & {\cdots} & 0 & 0 & {\cdots} & 0 \\
0 & 0 & 0 & {\alpha} & {\cdots} & 0 & 0 & {\cdots} & 0\\
{\vdots} & {\vdots} & {\vdots} & {\vdots} & {\ddots} & {\vdots} & {\vdots} & {\cdots} & {\vdots}\\
0 & 0 & 0 & 0 & {\cdots} & {\alpha} & 0 & {\cdots} & 0\\
0 & 0 & 0 & 0 & {\cdots} & 0 & {\frac{{\alpha}^2+2}{\alpha}} & {\cdots} & 0 \\
{\vdots} & {\vdots} & {\vdots} & {\vdots} & {\vdots} & {\vdots} & {\vdots} & {\ddots} & {\vdots}\\
0 & 0 & 0 & 0 & {\cdots} & 0 & 0 & {\cdots} & {\frac{{\alpha}^2+2}{\alpha}}
\end{bmatrix}
\end{equation*}

In this case also the formula $h={\lambda}+{\mu}$ and the notion of trace $h$ of the shape operator $S$ gives
\begin{equation*}
\begin{split}
h=&\alpha + \frac{{\alpha}^2+2}{\alpha}\\
=&(m+1){\alpha}+(m-2)\frac{{\alpha}^2+2}{\alpha}.
\end{split}
\end{equation*}
Then it implies that $(2m-3){\alpha}^2=-2m+4$, which gives us a contradiction for $m{\ge}3$.
\par
\vskip 6pt
Next we consider another case that the trace becomes $h=(m-1){\alpha}+(m-2)\frac{{\alpha}^2+2}{\alpha}$ in the above expression. Then $h={\lambda}+{\mu}={\alpha}+\frac{{\alpha}^2+2}{\alpha}$ gives
$$(2m-5){\alpha}^2+2(m-3)=0,$$
which also implies a contradiction for $m{\ge}4$.
\par
\vskip 6pt
Summing up the above discussions, we assert the following
\par
\vskip 6pt
\begin{MT1}\label{Theorem 1}\quad Let $M$ be a real hypersurface in complex quadric $Q^m$, $m{\ge}4$, with commuting Ricci tensor and $\mathfrak A$-isotropic normal. If the shape operator commutes with the structure tensor on the distribution ${\mathcal Q}^{\bot}$, then $M$ is locally congruent to a tube of radius $r$ over a totally geodesic ${\mathbb C}P^k$, $m=2k$, in $Q^{2k}$ or $M$ has $3$ distinct constant principal curvatures given by
$${\alpha}={\sqrt {2(m-3)}}, {\gamma}=0, {\lambda}=0,\  \text{and}\ \mu=-{\frac{2}{\sqrt {2(m-3)}}}\quad \text{or}$$
$${\alpha}={\sqrt {\frac{2}{3}(m-3)}}, {\gamma}=0, {\lambda}=0,\  \text{and}\ \mu=- \frac{\sqrt 6}{\sqrt {m-3}}$$
with corresponding principal curvature spaces
$$T_{\alpha}=[{\xi}], T_{\gamma}=[A{\xi}, AN], {\phi}(T_{\lambda})=T_{\mu}, \text{dim}\ T_{\lambda}=\text{dim}\ T_{\mu}=m-2.$$
\end{MT1}

\par
\vskip 6pt

\section{Proof of Main Theorem with $\frak A$-principal}\label{section 7}
\par
\vskip 6pt
In this section we want to prove our Main Theorem for real hypersurfaces with commuting Ricci tensor and $\frak A$-principal unit normal vector field.
\par
\vskip 6pt
From the basic formulas for the real structure $A$ and the K\"{a}hler structure $J$ we have the following
$$JAX=J\{BX+{\rho}(X)N\}={\phi}BX+{\eta}(BX)N-{\rho}(X){\xi},$$
$$AJX=A\{{\phi}X+{\eta}(X)N\}=B{\phi}X+{\rho}({\phi}X)N-{\eta}(X){\phi}B{\xi}-{\eta}(X){\eta}(B{\xi})N.$$
From this, the anti-commuting structure $AJ=-JA$ gives the following
\begin{equation}\label{e61}
{\phi}BX+{\eta}(BX)N-{\rho}(X){\xi}=-B{\phi}X-{\rho}({\phi}X)N+{\eta}(X){\phi}B{\xi}+{\eta}(X){\eta}(B{\xi})N.
\end{equation}
Then comparing the tangential and normal component of (\ref{e61}) gives the following respectively
\begin{equation}\label{e62}
{\eta}(BX)=-{\rho}({\phi}X)+{\eta}(X){\eta}(B{\xi}),
\end{equation}
and
\begin{equation}\label{e63}
{\phi}BX=-B{\phi}X+{\rho}(X){\xi}+{\eta}(X){\phi}B{\xi}.
\end{equation}

Since $N$ is $\mathfrak A$-principal, that is, $AN=N$, we know that $B{\xi}=-{\xi}$, and ${\phi}BX=-B{\phi}X$.  Then (\ref{e63}) gives
\begin{equation}\label{e64}
-2{\phi}BX+(Tr S)({\phi}S-S{\phi})X-({\phi}S^2-S^2{\phi})X=0,
\end{equation}
where we have used ${\rho}(X)=g(AX,N)=0$.  When $N$ is $\frak A$-principal, on the distribution ${\mathcal C}={\mathcal Q}$ we have

\begin{equation}\label{e65}
2S{\phi}S-{\alpha}({\phi}S+S{\phi})=2{\phi}.
\end{equation}
So if we put $SX={\lambda}X$ in (\ref{e65}), we have

\begin{equation}\label{e66}
S{\phi}X={\mu}{\phi}X=\frac{{\alpha}{\lambda}+2}{2{\lambda}-{\alpha}}{\phi}X .
\end{equation}
Then from (\ref{e63}) and (\ref{e65}) it follows that

\begin{equation}\label{e67}
-2{\phi}BX+ ({\lambda}-{\mu})\{h-({\lambda}+{\mu})\}{\phi}X=0.
\end{equation}

It is well known that the tangent space $T_{z}Q^m$ of the complex quadric $Q^m$ is decomposed as
$$T_{z}Q^m= V(A){\oplus}JV(A),$$
where $V(A)=\{X{\in}T_zQ^m{\vert}AX=X\}$ and $JV(A)=\{X{\in}T_zQ^m{\vert}AX=-X\}$. So $SX={\lambda}X$ for $X{\in}{\mathcal C}$
the vector field $X$ can be decomposed as follows:
$$X=Y+Z,\quad Y{\in}V(A), \quad Z{\in}JV(A),$$
where $AY=BY=Y$ and $AZ=BZ=-Z$. So it follows that $BX=AX=AY+AZ=Y-Z$. Then ${\phi}BX={\phi}Y-{\phi}Z$. From this, together with (\ref{e67}), it follows that
$$
-2({\phi}Y-{\phi}Z)+({\lambda}-{\mu})\{h-({\lambda}+{\mu})\}({\phi}Y+{\phi}Z)=0.$$
Then by taking inner products with ${\phi}Y$ and ${\phi}Z$ respectively, we get

\begin{equation}\label{e68}
({\lambda}-{\mu})\{h-({\lambda}+{\mu})\}-2=0
\end{equation}
and
\begin{equation}\label{e69}
({\lambda}-{\mu})\{h-({\lambda}+{\mu})\}+2=0.
\end{equation}
This gives a contradiction. Accordingly, we conclude that real hypersurfaces in $Q^m$ with commuting Ricci tensor and $\mathfrak A$-principal normal do not exist.
\par
\vskip 6pt
Summing up the above discussions, we assert the following
\par
\vskip 6pt
\begin{MT2}\label{Theorem 2}\quad There do not exist any real hypersurface in complex quadric $Q^m$, $m{\ge}4$, with commuting Ricci tensor and $\mathfrak A$-prinicipal normal vector field.
\end{MT2}
\par
\vskip 6pt
From Theorems 1 and 2, together with Lemma \ref{Lemma 4.2}, we give a complete proof of our Main Theorem in the introduction.
\vskip 6pt
\par

\begin{re}
 In this paper we have asserted that a tube over a totally geodesic ${\Bbb C}P^k$ in $Q^m$, $m=2k$, mentioned in our Main Theorem is Ricci commuting, that is, $Ric{\cdot}{\phi}={\phi}{\cdot}Ric$.
 But related to the notion of Ricci parallel, in the paper \cite{S5} we asserted that a tube over ${\Bbb C}P^k$ never has a parallel Ricci tensor, that is, the Ricci tensor does not satisfy ${\nabla}S=0$.
 \end{re}
\par

\end{document}